\documentstyle[preprint]{article}
\input{epsf}

\newtheorem{The}{Theorem}
\newtheorem{Lem}[The]{Lemma}
\newtheorem{Cor}[The]{Corollary}
\theorembodyfont{\rm}

\newtheorem{Def}{Definition}

\title {The Peano curve and counting occurrences of some patterns}
\author {Sergey Kitaev and Toufik Mansour \footnote{Research
financed by EC's
IHRP Programme, within the Research Training Network "Algebraic
Combinatorics in Europe", grant HPRN-CT-2001-00272}}
\address{Matematik, Chalmers tekniska h\"ogskola och G\"oteborgs
universitet, 412~96  G\"oteborg, Sweden}
\email{kitaev@math.chalmers.se, toufik@math.chalmers.se}

\def\mn{\mbox{-}}
\def\SS{{\mathcal S}}

\abstract{We introduce {\em Peano words}, which are words corresponding to finite approximations of the Peano space filling curve. We then find the number of occurrences of certain patterns in these words.}


\begin{document}
\maketitle
%-----------------------------------

\section{Introduction and Background}

We write permutations as words $\pi=a_1 a_2\cdots a_n$, whose letters are distinct and usually consist of the integers $1,2,\ldots,n$.

An occurrence of a pattern $p$ in a permutation $\pi$ is ``classically'' defined as a subsequence in $\pi$ (of the same length as the length of $p$) whose letters are in the same relative order as those in $p$. Formally speaking, for $r \leq n$, we say that a permutation $\sigma$ in the symmetric group ${\mathcal S}_n$ has an occurrence of the pattern $p \in {\mathcal S}_r$ if there exist $1 \leq i_1 < i_2 < \cdots < i_r \leq n$ such that $p = \sigma(i_1)\sigma(i_2) \ldots \sigma(i_r)$ in reduced form. The {\em reduced form} of a permutation $\sigma$ on a set $\{ j_1,j_2, \ldots ,j_r \}$, where $j_1 < j_2 < \cdots <j_r$, is a permutation ${\sigma}_1$ obtained by renaming the letters of the permutation $\sigma$ so that $j_i$ is renamed $i$ for all $i \in \{1, \ldots ,r\}$. For example, the reduced form of the permutation 3651 is 2431. The first case of classical patterns studied was that of permutations avoiding a pattern of length 3 in $\SS_3$. Knuth \cite{Knuth} found that, for any $\tau\in\SS_3$, the number $|\SS_n(\tau)|$ of $n$-permutations avoiding $\tau$ is $C_n$, the $n$th Catalan number.
Later, Simion and Schmidt \cite{SimSch} determined the number $|\SS_n(P)|$ of permutations in $\SS_n$ simultaneously avoiding any given set of patterns $P\subseteq\SS_3$.

In \cite{BabStein} Babson and Steingr\'{\i}msson introduced {\em generalised permutation patterns} that allow the requirement that two adjacent letters in a pattern must be adjacent in the permutation. In order to avoid confusion we write a "classical" pattern, say $231$, as $2$-$3$-$1$, and if we write, say $2$-$31$, then we mean that if this pattern occurs in the permutation, then the letters in the permutation that correspond to $3$ and $1$ are adjacent. For example, the permutation $\pi=516423$ has only one occurrence of the pattern $2$-$31$, namely the subword 564, whereas the pattern $2$-$3$-$1$ occurs, in addition, in the subwords 562 and 563. A motivation for introducing these patterns in \cite{BabStein} was the study of Mahonian statistics. A number of interesting results on generalised patterns were obtained in \cite{Claes}. Relations to several well studied combinatorial structures, such as set partitions, Dyck paths, Motzkin paths and involutions, were shown there.

Burstein \cite{Burstein} considered words instead of permutations. In particular, he found the number $|[k]^n(P)|$ of words of length $n$ in a $k$-letter alphabet that avoid all patterns from a set $P\subseteq\SS_3$ simultaneously. Burstein and Mansour \cite{BurMans1} (resp. \cite{BurMans2,BurMans3}) considered forbidden patterns (resp. generalized patterns) with repeated letters.

The most attention, in the papers on classical or generalized patterns, is paid to finding exact formulas and/or generating functions for the number of words or permutations avoiding, or having $k$ occurrences of, certain patterns. In~\cite{KitMans} the present authors suggested another problem, namely counting the number of occurrences of certain patterns in certain words. These words were chosen to be the set of all finite approximations of a sequence generated by a {\em morphism} with certain restrictions. A motivation for this choice was the interest in studying classes of sequences and words that are defined by iterative schemes~\cite{Lothaire,Salomaa}.

In the present paper we also study the number of occurrences of certain patterns in certain words. But here we choose these words to be the subdivision stages from which the {\em Peano curve} is obtained. We call these words the {\em Peano words}. The Peano curve was studied by the Italian mathematician Giuseppe Peano in 1890 as an example of a continuous space filling curve. We consider the Peano words and find the number of occurrences of the patterns $12$, $21$, $1^{\ell}$, $[x\mn y^{\ell})$, $(x^{\ell}\mn y]$ and $[x\mn y^{\ell}\mn z]$, where $x,y,z \in \{1,2,3\}$, $y^{\ell}=y\mn y\mn \cdots \mn y$ ($\ell$ times), and ``[`` in $p=[x-w)$ indicates that in an occurrence of $p$, the letter corresponding to the $x$ must be the first letter of the word.
%-----------------------------------

\section{The Peano curve and the Peano words}

We follow~\cite{GelbOlm} and present a description (of a curve that fills the unit square $S=[0,1]\times [0,1]$) given in 1891 by the German mathematician D. Hilbert.

As indicated in Figure~\ref{graph}, the idea is to subdivide $S$ and the unit interval $I=[0,1]$ into $4^n$ closed subsquares and subintervals, respectively, and to set up a correspondence between subsquares and subintervals so that inclusion relationships are preserved (at each stage of subdivision, if a square corresponds to an interval, then its subsquares correspond to subintervals of that interval).

We now define the continuous mapping $f$ of $I$ onto $S$: If $x\in
I$, then at each stage of subdivision $x$ belongs to {\em at
least} one closed subinterval. Select either one (if there are
two) and associate the corresponding square. In this way a
decreasing sequence of closed squares is obtained corresponding to
a decreasing sequence of closed intervals. This sequence of closed
squares has the property that there is exactly one point belonging
to all of them. This point is by definition $f(x)$. It can be
shown that the point $f(x)$ is well-defined, that is, independent
of any choice of intervals containing $x$; the range of $f$ is
$S$; and $f$ is continuous.
\begin{center}
\begin{figure}[h]
\begin{center}
\epsfxsize=5.9in \epsffile{peano.eps}\caption{the Peano words}
\label{graph}
\end{center}
\end{figure}
\end{center}
We now consider a subdivision stage (an iteration), go through the
curve inside $S$ starting in the point 1 (see Figure~\ref{graph}), and coding
any movement ``up'' by 1, ``right'' by 2, ''down'' by 3, ''left''
by 4. Thus, we start with the first iteration $X_1=123$, the
second iteration is $X_2=214112321233432$. More generally, it is
easy to see that the $n$-th iteration is given by
$$X_n=\varphi_1(X_{n-1})1X_{n-1}2X_{n-1}3\varphi_2(X_{n-1}),$$
where the function $\varphi_1(A)$ reverses the letters in the word $A$ and makes the substitution corresponding to the permutation $4123$, that is, 1 becomes 4 etc. The function $\varphi_2$ does the same, except with $4123$ replaced by $2341$. In this paper, we are interested in the words $X_n$, for $n=1,2,\ldots$, which appear as the subdivision stages of the Peano curve. We call these words the Peano words.

%-----------------------------------

\section{The main results}

It is easy to see that the length of the curve after the $n$-th iteration is $|X_n|=4^n-1$. Moreover, the following lemma holds.

\begin{Lem}\label{frequencies} The number of occurrences of the letters $1$, $2$, $3$ and $4$ in $X_n$ is given by $4^{n-1}$, $4^{n-1}+2^{n-1}-1$, $4^{n-1}$ and $4^{n-1}-2^{n-1}$ respectively. \end{Lem}
\begin{proof} Suppose $d^n_1$ (resp. $d^n_2$, $d^n_3$, $d^n_4$) denote the number of occurrences of the letter 1 (resp. 2,3,4) in the word $X_n$. It is easy to see, using the way we construct $X_n$, that
$$
\left(
\begin{array}{c}
d^n_1 \\
d^n_2 \\
d^n_3 \\
d^n_4
\end{array}
\right)=
\left(
\begin{array}{cccc}
2 & 1 & 0 & 1 \\
1 & 2 & 1 & 0 \\
0 & 1 & 2 & 1 \\
1 & 0 & 1 & 2
\end{array}
\right)
\left(
\begin{array}{c}
d^{n-1}_1 \\
d^{n-1}_2 \\
d^{n-1}_3 \\
d^{n-1}_4
\end{array}
\right)+
\left(
\begin{array}{c}
1 \\
1 \\
1 \\
0
\end{array}
\right).
$$
Using the diagonalization of the matrix in the identity above, namely the fact that
$$
\left(
\begin{array}{cccc}
2 & 1 & 0 & 1 \\
1 & 2 & 1 & 0 \\
0 & 1 & 2 & 1 \\
1 & 0 & 1 & 2
\end{array}
\right)=
\left(
\begin{array}{rrrr}
-1 & -1 & 0 & 1 \\
1 & 0 & -1 & 1 \\
-1 & 1 & 0 & 1 \\
1 & 0 & 1 & 1
\end{array}
\right)
\left(
\begin{array}{cccc}
0 & 0 & 0 & 0 \\
0 & 2 & 0 & 0 \\
0 & 0 & 2 & 0 \\
0 & 0 & 0 & 4
\end{array}
\right)
\left(
\begin{array}{cccc}
-1/4 & 1/4 & -1/4 & 1/4 \\
-1/2 & 0 & 1/2 & 0 \\
0 &-1/2 & 0 & 1/2 \\
1/4 & 1/4 & 1/4 & 1/4
\end{array}
\right),
$$
we get that the vector $(d^n_1, d^n_2, d^n_3, d^n_4)^{'}$ is equal to
$$
\left(
\begin{array}{rrrr}
-1 & -1 & 0 & 1 \\
1 & 0 & -1 & 1 \\
-1 & 1 & 0 & 1 \\
1 & 0 & 1 & 1
\end{array}
\right)
\left(
\begin{array}{cccc}
1 & 0 & 0 & 0 \\
0 & 2^n-1 & 0 & 0 \\
0 & 0 & 2^n-1 & 0 \\
0 & 0 & 0 & (4^n-1)/3
\end{array}
\right)
\left(
\begin{array}{cccc}
-1/4 & 1/4 & -1/4 & 1/4 \\
-1/2 & 0 & 1/2 & 0 \\
0 &-1/2 & 0 & 1/2 \\
1/4 & 1/4 & 1/4 & 1/4
\end{array}
\right)
\left(
\begin{array}{c}
1 \\
1 \\
1 \\
0
\end{array}
\right),
$$
which is equal to the vector $(4^{n-1}, 4^{n-1}+2^{n-1}-1, 4^{n-1}, 4^{n-1}-2^{n-1})^{'}$.
\end{proof}

As a corollary to Lemma~\ref{frequencies} we have the following.

\begin{Cor}\label{cor2}
The number of occurrences of the pattern $\tau=\underbrace{1\mn 1\mn \cdots \mn 1}_{\ell}=1^{\ell}$ in $X_n$ is equal to
$$
{4^{n-1}-2^{n-1} \choose \ell}+2{4^{n-1}\choose \ell}+{4^{n-1}+2^{n-1}-1\choose \ell}.
$$
\end{Cor}

\begin{proof}
The number of occurrences of a subsequence $\underbrace{i\mn i\mn \cdots \mn i}_{\ell}$ in $X_n$, for $i=1,2,3,4$, is obviously given by ${d^n_i \choose \ell}$, where $d^n_i$ is defined and determined in the proof of Lemma~\ref{frequencies}. The rest is easy to see.
\end{proof}

\begin{Def} Let $r(A)$ (resp. $d(A)$) denote the number of occurrences of the pattern $12$ (resp. $21$), that is the number of rises (resp. descents), in a word $A$. \end{Def}

\begin{Lem}\label{aux1} Suppose $A=1X3$ and $B=2Y2$ for some words $X$ and $Y$. Then $r(\varphi_1(A))=d(A)+1$, $d(\varphi_1(A))=r(A)-1$, $r(\varphi_2(B))=d(B)$ and $d(\varphi_2(B))=r(B)$.
\end{Lem}

\begin{proof}
If $\bar{A}$ and $\bar{B}$ denote the reverses of $A$ and $B$ then $r(\bar{A})=d(A)$, $d(\bar{A})=r(A)$, $r(\bar{B})=d(B)$, and $d(\bar{B})=r(B)$.

We consider two factorizations of each word $\bar{A}$ and $\bar{B}$. We can write $\bar{A}$ as
$$\bar{A}=3A_1\underbrace{1\ldots 1}_{i_1}A_2\underbrace{1\ldots 1}_{i_2}A_3\ldots A_{k}\underbrace{1\ldots 1}_{i_k},$$
where $A_i$, for $i=1,2,\ldots, k$ is a word over the alphabet $\{ 2,3,4 \}$, only $A_1$ can be the empty word $\epsilon$, and $i_j\geq 1$ for $j=1,2\ldots,k$. Also, we can write $\bar{A}$ as
$$\bar{A}=3A^{'}_0\underbrace{4\ldots 4}_{i^{'}_1}A^{'}_1\underbrace{4\ldots 4}_{i^{'}_2}A^{'}_2\ldots A^{'}_{k-1}\underbrace{4\ldots 4}_{i^{'}_k}A^{'}_{k}1,$$
where $A^{'}_i$, for $i=0,1,\ldots, k$ is a word over the alphabet $\{ 1,2,3 \}$, only $A^{'}_0$ and $A^{'}_k$ can be $\epsilon$, and $i^{'}_j\geq 1$ for $j=1,2\ldots,k$.

The word $\bar{B}$ can be written as
$$\bar{B}=2B_0\underbrace{1\ldots 1}_{j_1}B_1\underbrace{1\ldots 1}_{j_2}B_2\ldots B_{\ell-1}\underbrace{1\ldots 1}_{j_{\ell}}B_{\ell}2,$$
where $B_i$, for $i=0,1,\ldots, \ell$, is a word over the alphabet $\{ 2,3,4 \}$, only $B_0$ and $B_{\ell}$ can be $\epsilon$, and $j_i\geq 1$ for $i=1,2\ldots,\ell$. Also, $\bar{B}$ can be written as

$$\bar{B}=2B^{'}_0\underbrace{4\ldots 4}_{j^{'}_1}B^{'}_1\underbrace{4\ldots 4}_{j^{'}_2}B^{'}_2\ldots B^{'}_{\ell-1}\underbrace{4\ldots 4}_{j^{'}_{\ell}}B^{'}_{\ell}2,$$ where $B^{'}_i$, for $i=0,1,\ldots, \ell$, is a word over the alphabet $\{ 1,2,3 \}$, only $B^{'}_0$ and $B^{'}_{\ell}$ can be $\epsilon$, and $j^{'}_i\geq 1$ for $i=1,2\ldots,\ell$.

It follows from the definitions that  $\varphi_1(A)$ and $\varphi_1(B)$ (resp. $\varphi_2(A)$ and $\varphi_2(B)$) are obtained from $\bar{A}$ and $\bar{B}$ by permuting the letters with the function $\pi_1$ (resp. $\pi_2$) that acts as the permutation $4123$ (resp. $2341$).

We now consider the first factorizations of $\bar{A}$ and $\bar{B}$, and the function $\pi_1$. It is easy to see that if $W$ is equal to $A_i$, or $B_i$, or $3A_1$, or $2B_0$, or $B_{\ell}2$, then $r(W)=r(\pi_1(W))$ and $d(W)=d(\pi_1(W))$, since $\pi_1$ is an order-preserving function when it acts from the set $\{ 2,3,4 \}$ to the set $\{ 1,2,3 \}$. From the other hand, occurrences of the rises $12$, $13$ and $14$ (resp. the descents $41$, $31$ and $21$) in $\bar{A}$ and $\bar{B}$, give occurrences of the descents $41$, $42$ and $43$ (resp. the rises $34$, $24$ and $14$) in $\pi_1(\bar{A})$ and $\pi_1(\bar{B})$ respectively. If we now read the first factorizations of $\bar{A}$ and $\bar{B}$ from the left to the right, then the occurrences of the subwords $a1$ alternate with the occurrences of the subwords $1b$, where $a,b \in \{ 2,3,4 \}$. Moreover, in the factorization of $\bar{A}$, we begin and end with the subword $a1$ for some $a\in \{ 2,3,4 \}$, which gives that $d(A)+1=r(\bar{A})+1=r(\pi_1(\bar{A}))=r(\varphi_1(A))$ and $r(A)-1=d(\bar{A})-1=d(\pi_1(\bar{A}))=d(\varphi_1(A))$; in the factorization of $\bar{B}$, we begin with the subword $a1$ and end with the subword $1b$ for some $a,b\in \{ 2,3,4 \}$, which gives that $d(B)=r(\bar{B})=r(\pi_1(\bar{B}))=r(\varphi_1(B))$ and $r(B)=d(\bar{B})=d(\pi_1(\bar{B}))=d(\varphi_1(B))$.

If we consider the second factorizations of $\bar{A}$ and $\bar{B}$, and the function $\pi_2$, one can see that if $W$ is equal to $A^{'}_i$, or $B^{'}_i$, or $3A^{'}_0$, or $A^{'}_k1$, or $2B^{'}_0$, or $B^{'}_{\ell}2$, then $r(W)=r(\pi_2(W))$ and $d(W)=d(\pi_2(W))$, since $\pi_2$ is an order-preserving function when it acts from the set $\{ 1,2,3 \}$ to the set $\{ 2,3,4 \}$. From the other hand, occurrences of the rises $14$, $24$ and $34$ (resp. the descents $41$, $42$ and $43$) in $\bar{A}$ and $\bar{B}$, give occurrences of the descents $21$, $31$ and $41$ (resp. the rise $12$, $13$ and $14$) in $\pi_2(\bar{A})$ and $\pi_2(\bar{B})$ respectively. If we now read the second factorizations of $\bar{A}$ and $\bar{B}$ from the left to the right, then the occurrences of the subwords $a4$ alternate with the occurrences of the subwords $4b$, where $a,b \in \{ 1,2,3 \}$. Moreover, in both cases, we begin with the subword $a4$ and end with the subword $4b$ for some $a,b\in \{ 1,2,3 \}$, which gives that $d(A)=r(\bar{A})=r(\pi_2(\bar{A}))=r(\varphi_2(A))$, $r(A)=d(\bar{A})=d(\pi_2(\bar{A}))=d(\varphi_2(A))$, $d(B)=r(\bar{B})=r(\pi_2(\bar{B}))=r(\varphi_2(B))$ and $r(B)=d(\bar{B})=d(\pi_2(\bar{B}))=d(\varphi_2(B))$.
\end{proof}

\begin{The}\label{theorem4} Let $r_n$ {\rm(}resp. $d_n${\rm)} be the number of occurrences of the pattern $12$ {\rm(}resp. $21${\rm)} in $X_n$. Then for all $k\geq0$,
$$
\begin{array}{l}
r_{2k+1}=\frac{2}{5}(4\cdot 16^k + 1), \\[2mm]
r_{2k+2}=\frac{2}{5}(16^{k+1}-1), \\[2mm]
d_{2k+1}=\frac{8}{5}(16^k-1), \\[2mm]
d_{2k+2}=\frac{2}{5}(16^{k+1}-1).
\end{array}
$$
\end{The}

\begin{proof} Using the properties of $\varphi_1$ and $\varphi_2$, as well as the way we construct $X_n$, it is easy to check by induction, that $X_{2k+1}$ and $X_{2k+2}$ can be factorized as follow:
$$X_{2k+1}=\underbrace{1X^{(1)}1}_{\varphi_1(X_{4k})}1\underbrace{2Y^{(1)}2}_{X_{4k}}2\underbrace{2Y^{(1)}2}_{X_{4k}}3\underbrace{3Z^{(1)}3}_{\varphi_2(X_{4k})},$$
$$X_{2k+2}=\underbrace{2X^{(2)}4}_{\varphi_1(X_{4k+1})}1\underbrace{1Y^{(2)}3}_{X_{4k+1}}2\underbrace{1Y^{(2)}3}_{X_{4k+1}}3\underbrace{4Z^{(2)}2}_{\varphi_2(X_{4k+1})},$$
where $X^{(i)}$, $Y^{(i)}$ and $Z^{(i)}$ are some words for $i=1,2$.

Suppose we know $r_{2k+1}$ and $d_{2k+1}$ for some $k$. Since $X_{2k+1}=1W3$ for some word $W$, using Lemma~\ref{aux1} and the factorization of the word $X_{2k+2}$, we can find $r_{2k+2}$ and $d_{2k+2}$. Indeed, $\varphi_1(X_{4k+1})$ has $d_{2k+1}+1$ rises and $r_{2k+1}-1$ descents; $\varphi_2(X_{4k+1})$ has $d_{2k+1}$ rises and $r_{2k+1}$ descents; two subwords $X_{2k+1}$ give $2r_{2k+1}$ rises and $2d_{2k+1}$ descents. Besides, we have some extra rises and descents appeared between different blocks of the decomposition. They are one extra rise between the letter 3 and the subword $\varphi_2(X_{4k+1})$, and 3 extra descents between the subword $\varphi_1(X_{4k+1})$ and the letter 1, the subword $X_{4k+1}$ and the letter 2, the letter 2 and the subword $X_{4k+1}$. Thus, $r_{2k+2}=2r_{2k+1}+2d_{2k+1}+2$ and $d_{2k+2}=2r_{2k+1}+2d_{2k+1}+2$, which shows, in particular, that for even $n$, in $X_n$, the number of rises is equal to the number of descents.

We now analyze the factorization of $X_{2k+3}$, which is similar to that of $X_{2k+1}$. Using the fact that $X_{2k+2}=2W^{'}2$ for some word $W^{'}$ and Lemma~\ref{aux1}, we can find $r_{2k+3}$ and $d_{2k+3}$. Indeed, we can use the similar considerations as above to get $r_{2k+3}=2r_{2k+2}+2d_{2k+2}+2=8r_{2k+1}+8d_{2k+1}+10$ and $d_{2k+3}=2r_{2k+2}+2d_{2k+2}=8r_{2k+1}+8d_{2k+1}+8$. Thus, if $x_{k}$ denote the vector $(r_{2k+1},d_{2k+1})^{'}$ then
$$
x_{k+1}=
\left(
\begin{array}{cc}
8 & 8 \\
8 & 8 \\
\end{array}
\right)
x_k+
\left(
\begin{array}{c}
10 \\
8 \\
\end{array}
\right),
$$
with $x_0=(2,0)$, since in $X_1=123$, there are two rises and no descents. This recurrence relation, using diagonalization of the matrix in it, leads us to $$x_k=(\frac{2}{5}(4\cdot 16^k+1), \frac{8}{5}(16^k-1))^{'}.$$

Finally, $r_{2k+2}=d_{2k+2}=2r_{2k+1}+2d_{2k+1}+2=\frac{2}{5}(16^{k+1}-1)$.
\end{proof}

Let $N_{\tau}(W)$ denote the number of occurrences of the pattern $\tau$ in the word $W$.

Using Lemma~\ref{frequencies} and the proof of Theorem~\ref{theorem4}, we can count, for $X_n$, the number of occurrences of the patterns $\tau_1(x,y)=[x\mn y^{\ell})$, $\tau_2(x,y)=x^{\ell}\mn y]$ and $\tau_3(x,y,z)=[x\mn y^{\ell} \mn z]$, where $x,y,z \in \{1,2,3\}$, $y^{\ell}=y\mn y\mn \cdots \mn y$ ($\ell$ times), and ``[`` in $p=[x-w)$ indicates that in an occurrence of $p$, the letter corresponding to the $x$ must be the first letter of the word, whereas ``]'' in $\tau_3(x,y,z)$ indicates that in an occurrence of $\tau_3(x,y,z)$, the letter corresponding to the $z$ must be the last (rightmost) letter of the word.

If we consider, for instance, the pattern $\tau_1(1,2)=[1\mn 2^{\ell})$ then the letter 1 in this pattern must correspond to the leftmost letter of the word $X_n$. Now if $n=2k+1$ then from the proof of Theorem~\ref{theorem4} $X_n=1W$ for some word $W$, which means that to the sequence $2^{\ell}$ there can correspond any subsequence $i^{\ell}$ in $X_n$, where $i=2,3,4$. Thus, using Lemma~\ref{frequencies} and the way we prove Corollary~\ref{cor2}, there are ${4^{2k}-2^{2k} \choose \ell}+{4^{2k}\choose \ell}+{4^{2k}+2^{2k}-1\choose \ell}$ occurrences of the pattern $\tau_1(1,2)$ in $X_{2k+1}$. If $n=2k+2$ then $X_n=2W$ for some word $W$ and for the sequence $2^{\ell}$ there correspond any subsequence $i^{\ell}$ in $X_n$, where $i=3,4$. Thus, $N_{\tau_1(1,2)}(X_{2k+2})={4^{2k}-2^{2k} \choose \ell}+{4^{2k}\choose \ell}$.

In the example above, as well as in the following considerations, we assume $\ell$ to be greater then $0$. If $\ell=0$ then obviously $N_{\tau_1(x,y)}(X_n)=N_{\tau_2(x,y)}(X_n)=1$, whereas $N_{\tau_3(x,y,z)}(X_n)$ is equal to $1$ if $x<z$ and $n=2k+1$, or $x=z$ and $n=2k+2$, and it is equal to 0 otherwise.

When we consider $\tau_3(x,y,z)(X_n)$, we observe that since $X_{2k+2}=2W2$ for some $W$, $N_{\tau_3(x,y,z)}(X_{2k+2})=0$, whenever $x\neq z$. Also, since $X_{2k+1}=1W3$ for some $W$, $N_{\tau_3(x,y,z)}(X_{2k+1})=0$, whenever $x\geq z$.

Let us consider the pattern $\tau_3(2,1,3)=[2\mn 1^{\ell} \mn 3]$. As it was mentioned before, $N_{\tau_3(2,1,3)}(X_{2k+2})=0$. But, if we consider $X_{2k+1}=1W3$, then it is easy to see that $N_{\tau_3(2,1,3)}(X_{2k+1})=0$, since the leftmost letter of $X_{2k+1}$ is the least letter, which means that it cannot correspond to the letter $2$ in the pattern. As one more example, we can consider the pattern $\tau_3(1,1,2)=[1\mn 1^{\ell} \mn 2]$. We are only interested in case $X_n=X_{2k+1}$, since $N_{\tau_3(1,1,2)}(X_{2k+2})=0$. The number of occurrences of the pattern is obviously given by the number of ways to choose $\ell$ letters among $4^{2k}-1$ letters $1$ (totally, there are $4^{2k}$ letters 1 according to Lemma~\ref{frequencies}, but we cannot consider the leftmost 1 since it corresponds to the leftmost 1 in the pattern). Thus, $N_{\tau_3(1,1,2)}(X_{2k+1})={4^{2k}-1 \choose \ell}$.

All the other cases of $x$, $y$, $z$ in the patterns $\tau_1(x,y)$, $\tau_2(x,y)$ and $\tau_3(x,y,z)$ can be considered in the same way. Let $S_1$ and $S_2$ denote the following:
$$S_1={4^{2k}-2^{2k} \choose \ell}+ {4^{2k} \choose \ell}+{4^{2k}+2^{2k}-1 \choose \ell},\mbox{\ \ }S_2={4^{2k+1} \choose \ell}+{4^{2k+1}-2^{2k+1} \choose \ell}.$$
The tables below give all the results concerning the patterns under consideration, except those triples $(x,y,z)$, for which $N_{\tau_3(x,y,z)}(X_{n})=0$ for all $n$.

\begin{center}
\begin{tabular}{|c|c|c|c|c|c|}
\hline
$x$ & $y$ & $N_{\tau_1(x,y)}(X_{2k+1})$ & $N_{\tau_2(x,y)}(X_{2k+1})$ & $N_{\tau_1(x,y)}(X_{2k+2})$ & $N_{\tau_2(x,y)}(X_{2k+2})$ \\[1,5mm]
\hline
1 & 1 & ${4^{2k}-1 \choose \ell}$ & ${4^{2k}-1 \choose \ell}$ & ${4^{2k+1}+2^{2k+1}-1 \choose \ell}$ & ${4^{2k+1}+2^{2k+1}-1 \choose \ell}$ \\[1,5mm]
\hline
1 & 2 & $S_1$ & ${4^{2k} \choose \ell}+{4^{2k}+2^{2k}-1 \choose \ell}$ & $S_2$ & ${4^{2k+1} \choose \ell}$ \\[1,5mm]
\hline
2 & 1 & 0 & ${4^{2k} - 2^{2k}\choose \ell}$ & ${4^{2k+1} \choose \ell}$ & $S_2$ \\[1,5mm]
\hline
\end{tabular}

\

\

\begin{tabular}{|c|c|c|c|c|}
\hline
$x$ & $y$ & z & $N_{\tau_3(x,y,z)}(X_{2k+1})$ & $N_{\tau_3(x,y,z)}(X_{2k+2})$ \\[1,5mm]
\hline
1 & 1 & 1 & 0 & ${4^{2k+1}-2 \choose \ell}$ \\[1,5mm]
\hline
1 & 1 & 2 & ${4^{2k}-1 \choose \ell}$ & 0 \\[1,5mm]
\hline
1 & 2 & 1 & 0 & $S_2$ \\[1,5mm]
\hline
1 & 2 & 2 &${4^{2k}-1 \choose \ell}$ & 0 \\[1,5mm]
\hline
2 & 1 & 2 & 0 & ${4^{2k+1} \choose \ell}$ \\[1,5mm]
\hline
1 & 2 & 3 & ${4^{2k}+2^{2k}-1 \choose \ell}$ & 0 \\[1,5mm]
\hline
1 & 3 & 2 & ${4^{2k}-2^{2k} \choose \ell}$ & 0 \\[1,5mm]
\hline
\end{tabular}
\end{center}

%-----------------------------------

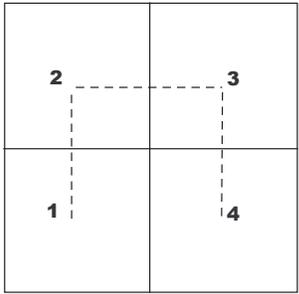
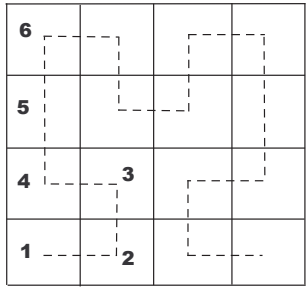
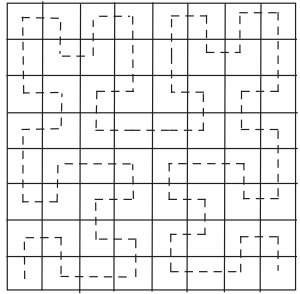
\end{document}